\documentclass[12pt,reqno]{amsart}

\usepackage{dsfont, amsfonts, amsmath, amssymb,amscd, stmaryrd, latexsym, dsfont,amssymb,amsfonts,hyperref,enumerate,amscd,tikz}

\newtheorem{thm}{Théorème}[section]

\newtheorem{cor}[thm]{Corollaire}
\newtheorem{lem}[thm]{Lemme}
\newtheorem{prop}[thm]{Proposition}
\newtheorem{defn}[thm]{Définition}
\newtheorem{rem}[thm]{Remarque}
\newtheorem{rems}[thm]{Remarques}
\newtheorem{ex}[thm]{Exemple}
\newtheorem{exs}[thm]{Exemples}

\usepackage{geometry}
\begin{document}

  \title[Survey on generalizations of Hopficity of modules]{Survey on generalizations of Hopficity of modules}
  \author[Abderrahim El moussaouy]{Abderrahim El moussaouy}
  \address{Abderrahim El moussaouy: Department of Mathematics, Faculty of Sciences,
  	University of Mohammed First, Oujda, Morocco }
  \email{a.elmoussaouy@ump.ac.ma
  }


\maketitle
\begin{abstract}
	
	The main aim of this paper is the  Hopficity of module classes, the study of modules (rings) by properties of their endomorphisms is a
	classical research subject. In 1986, Hiremath \cite{Hi} introduced the concepts of Hopfian
	modules and rings, the notion of Hopfian modules are defined as a generalization of modules of
	finite length as the modules whose surjective endomorphisms are isomorphisms. Later, the dual concepts co-Hopfian modules and rings were
	given. Hopfian and co-Hopfian modules (rings) have been investigated by several authors. For example, Hiremath \cite{Hi}, Varadarajan \cite{Va}, \cite{Va1}, Xue \cite{Xu}, Haghany \cite{Hag},
	Liu \cite{Li}, and Yang and Liu \cite{Yl}. In 2001, Haghany and Vedadi \cite{Ha}, and in 2002, Ghorbani
	and Haghany \cite{Gh}, respectively, introduced and investigated the weakly co-Hopfian
	and generalized Hopfian modules. These modules and several generalizations of them are extensively studied also by several authors.

\end{abstract}

\newpage

\section*{Résumé}	

L'objectif principal de cet article est d'étudier les classes de Hopficité de modules. L'étude des modules (anneaux) par les propriétés de leurs endomorphismes est un sujet de recherche classique. En 1986, Hiremath \cite{Hi} a introduit les concepts de modules Hopfiens et d'anneaux Hopfiens, la notion de modules Hopfiens est définie comme une généralisation de modules de longueur finie comme les modules dont les endomorphismes surjectifs sont des isomorphismes. Plus tard, les deux concepts de modules co-Hopfiens et d'anneaux co-Hopfiens ont été donnés. Les modules (anneaux) Hopfiens et co-Hopfiens ont été étudiés par plusieurs auteurs. Par exemple, Hiremath \cite{Hi}, Varadarajan \cite{Va}, \cite{Va1}, Xue \cite{Xu}, Haghany \cite{Hag}, Liu \cite{Li} et Yang et Liu \cite{Yl}. En 2001, Haghany et Vedadi, \cite{Ha}, et en 2002, Ghorbani et Haghany \cite{Gh}, respectivement, ont introduit et étudié les modules faiblement co-Hopfiens et Hopfiens généralisés. Ces modules et plusieurs généralisations ont été étudiés également par plusieurs auteurs. \\


\newpage

\section*{Notations}
Nous fixons les notations suivantes. Soient $A$ un anneau associatif
unitaire, $M$ un $A$-module \`{a} droite et $N$ un
sous-module de $M$.
{\begin{enumerate}[\rm $\bullet$]
		\item $\mathrm{End}(M)$ l'anneau des endomorphismes de $M$.
		\item $N\leq M$: $N$ est un sous-module de $M$,
		\item $N\unlhd M$: $N$ est un sous-module totalement invariant de $M$ $(f(N)\leq N$ pour tout $f\in \mathrm{End}(M))$,
		\item $N\leq^{\oplus} M$: $N$ est un facteur direct de $M$,
		\item $N\leq^{e} M$: $N$ est un sous-module essentiel de $M$,
		\item $N\ll M$: $N$ est un sous-module superflu dans $M$,
		\item $N\ll_{\mu} M$: $N$ est un sous-module $\mu$-superflu dans $M$,
		\item $N\ll_{\delta} M$: $N$ est un sous-module $\delta$-superflu dans $M$,
		\item $N\ll_{\gamma} M$: $N$ est un sous-module $\gamma$-superflu dans $M$,
		\item $N\ll_{J} M$: $N$ est un sous-module Jacobson-superflu dans $M$,
		\item $Z(M)$: le sous-module singulier de $M$,
		\item $Z^{*}(M)$: le sous-module cosingulier de $M$,
		
		\item $\mathrm{Rad}(M)$ le radical de Jacobson de $M$.
		\item $J(A)$ le radical de Jacobson de $A$.
		\item $\delta(A)$ l’intersection de tous les idéaux essentiels maximaux de $A$.
		\item $\mathrm{E}(M)$ l'enveloppe injective de $M$.
		\item CCA la condition de chaîne ascendante.
		\item CCD la condition de chaîne descendante.
		\item $A_{p}=\{\frac{a}{p^{n}}$; $a\in \mathds{Z}$ et $n \in \mathds{N}\}$, où $p$ est un nombre premier.
		\item $\mathbb{Z}_{p^{\infty}}=A_{p}/\mathds{Z}$.
		\item $M_{n}(A)$: l’ensemble des matrices carrées d’ordre $n$ à coefficients dans $A$.\\
		\\

		Dans cet article, tous les anneaux sont considérés
		associatifs unitaires et les modules sont considérés des modules à
		droite sauf mention du contraire.
\end{enumerate}}

\newpage

\section{Introduction}

Au début des années quatre-vingt, A. Kaidi et M. Sangharé ont introduit , la notion des modules vérifiant les
propriétés (I), (S) et (F), \cite{Amin}. On dit qu'un A-module à
droite $M$ vérifie la propriété (I) (resp. (S)) si tout
endomorphisme injectif (resp. surjectif) de $M$  est un automorphisme
de $M$, on dit que $M$ vérifie  la propriété (F) si pour
tout endomorphisme $f$ de $M$, il existe un entier $n\geq 1$ tel que
$M=Ker f^{n}\oplus Im f^{n}$.

En 1986, la notion du module vérifiant la condition (S) a été nommée par Hiremath "module Hopfien", \cite{Hi}. Un peu plus tard la notion du module vérifiant la condition (I) a été
nommée par Varadarajan "module co-Hopfien",\cite{Va1}.

On dit qu'un sous-module $N$ d'un $A$-module $M$ est essentiel de
$M$ ($N\leq^{e} M$) si $N\cap L=0$ implique $L=0$, pour tout
sous-module $L$ de $M$. En 2001, \cite{Ha}, A. Haghany et M. R.
Vedadi ont introduit la notion du module faiblement co-Hopfien. Un
$A$-module $M$ est appelé faiblement co-Hopfien, si pour tout
endomorphisme injectif $f$ de $M$, l'image de $f$ est essentiel de
$M$ $(Im f \leq^{e} M)$.

On dit qu'un sous-module $N$ d'un $A$-module $M$ est superflu dans
$M$ ($N\ll M$) si $N+L=M$ implique $L=M$, pour tout sous-module $L$
de $M$. En 2002, \cite{Gh}, A. Ghorbani et A. Haghany ont introduit
la notion du module Hopfien généralisé. On dit qu'un
$A$-module $M$ est Hopfien généralisé, si pour tout
endomorphisme surjectif $f$ de $M$, le noyau de $f$ est superflu
dans $M$ $(Ker f \ll M)$.

En 2005, \cite{Wa}, Y. Wang a introduit la notion du module
co-Hopfien généralisé et la notion du module faiblement
Hopfien. On dit qu'un $A$-module $M$ est faiblement Hopfien, si tout
endomorphisme surjectif superflu $f$ de $M$ est un automorphisme. Et
on dit qu'un $A$-module $M$ est co-Hopfien généralisé,
si tout endomorphisme injectif essentiel $f$ de $M$ est un
automorphisme.

En 2007, \cite{Hm}, A. Hmaimou, A. kaidi et E. Sanchez Campos ont
introduit la notion du module fortement Hopfien et la notion du module fortement
co-Hopfien. On dit qu'un $A$-module $M$ est fortement Hopfien, si
pour tout endomorphisme $f$ de $M$ la suite croissante: $Ker f
\subseteq Ker f^{2} \subseteq  ... \subseteq Ker f^{n} \subseteq
...$ est stationnaire. Et on dit qu'un $A$-module $M$ est fortement
co-Hopfien, si pour tout endomorphisme $f$ de $M$ la suite
décroissante: $Im f \supseteq Im f^{2} \supseteq  ... \supseteq
Im f^{n} \supseteq ...$ est stationnaire.

En 2008, \cite{Pi}, P. Aydogdu et  A.C. Ozcan ont introduit
la notion du module semi Hopfien et la notion du module semi co-Hopfien. On dit qu'un
$A$-module $M$ est semi Hopfien, si pour tout endomorphisme surjectif $f$ de $M$, le noyau de $f$ est un facteur direct de $M$. Et on dit qu'un
$A$-module $M$ est semi co-Hopfien, si pour tout endomorphisme injectif $f$ de $M$, l'image de $f$ est un facteur direct de $M$.

De tels modules et d'autres généralisations ont été
introduits et étudiés par plusieurs auteurs (\cite{Pi}, \cite{A}, \cite{AB}, \cite{ABC},
\cite{Ga}, \cite{Gha}, \cite{Ha}, \cite{Hi},\cite{Hm}, \cite{Va1}, \cite{Fer}, \cite{Alg}, \cite{Ali}, \cite{Kho}, \cite{Jor},\cite{Ico},\cite{Tur},
\cite{Wa}).

\section{Modules Hopfien et co-Hopfien}

\begin{defn}\cite{Hi}.
	Un $A$-module $M$ est dit Hopfien si tout endomorphisme surjectif de
	$M$ est bijectif.
\end{defn}
\begin{defn}\cite{Va1}.
	Un $A$-module $M$ est dit co-Hopfien si tout endomorphisme injectif
	de $M$ est bijectif.
\end{defn}
\begin{prop}\cite{Vas}.
	Si $A$ est un anneau commutatif, alors tout $A$-module de type fini
	est Hopfien.
\end{prop}

\begin{rem}\cite{Ri}.
	Tout module noethérien (resp, artinien) est Hopfien (resp,
	co-Hopfien).
\end{rem}
Un module Hopfien (resp, co-Hopfien) n'est pas en général
noethérien (resp, artinien) comme le montre l'exemple suivant:

\begin{ex}\cite{Hm}.
	Le groupe additif $\mathbb{Q}$ des nombres rationnels est Hopfien et
	co-Hopfien mais n'est ni noethérien ni artinien.
\end{ex}
\begin{defn}\cite{Br,Mo}.
	L'anneau $A$ est dit Dedekind Fini si pour tout $a,b \in A$, $ab=1
	\Rightarrow ba=1$. Le module $M$ est dit Dedekind Fini si l'anneau
	$\mathrm{End}(M)$ est Dedekind Fini. On peut vérifier que $M$ est
	Dedekind Fini si et seulement si $M$ n'est pas isomorphe à un
	facteur direct propre de lui-même.
\end{defn}
\begin{prop}\cite{Br}.
	Soit $M$ un $A$-module Hopfien ou co-Hopfien, alors $M$ est Dedekind
	Fini.
\end{prop}
La réciproque n'est pas en général vraie comme le montre
l'exemple suivant:

\begin{ex}
	Le groupe abélien $\mathbb{Z}$ est Dedekind Fini mais n'est pas
	co-Hopfien et le groupe $\mathbb{Z}_{p^{\infty}}$ est Dedekind Fini
	mais n'est pas Hopfien.
\end{ex}

\begin{defn}.
	Un $A$-module $E$ est dit injectif si pour tout
	homomorphisme injectif $g$ de M vers N  et pour tout homomorphisme, $\gamma$
	de M vers E, il existe un homomorphisme $h$ de N vers E tel que :
	$\gamma = hg$ (i.e., il existe $h:N \rightarrow E$ tel que le
	diagramme
	
	\begin{center}
		\begin{tikzpicture}
			\begin{scope}[xscale=2.5,yscale=2.5]
				\node (A) at (0,1)   {$M$};
				\node (B) at (1,1){$N$};
				\node(C) at (-1,1) {$0$};
				
				\node (E)   at (0,0) {$E$};
				\draw [->,>=latex] (C) -- (A)   node[midway,below,rotate=0]  {$ $};
				
				\draw [->,>=latex] (A) -- (B)   node[midway,above,rotate=0]  {$g$};
				\draw [->,>=latex] (A) -- (E)   node[midway,below,rotate=0]  {$ $ $ $ $\gamma$};

				\draw [->,>=latex,dotted] (B) -- (E)   node[midway,below,rotate=0]  {$ $ $ $ $h$};
				
			\end{scope}
		\end{tikzpicture}
	\end{center}

	est commutatif).
\end{defn}
\begin{defn}.
	Un $A$-module $P$ est dit projectif si pour tout
	homomorphisme surjectif $g$ de M vers N  et pour tout homomorphisme, $\gamma$
	de P vers N, il existe un homomorphisme $h$ de P vers M tel que :
	$\gamma = gh$ (i.e., il existe $h:P \rightarrow M$ tel que le
	diagramme
	
	\begin{center}
		\begin{tikzpicture}
			\begin{scope}[xscale=2.5,yscale=2.5]
				\node (A) at (1,0)   {$N$};
				\node (B) at (0,0){$M$};
				\node(C) at (2,0) {$0$};
				\node (E)   at (1,1) {$P$};
				\draw [->,>=latex] (A) -- (C)   node[midway,below,rotate=0]  {$ $};
				\draw [->,>=latex] (B) -- (A)   node[midway,above,rotate=0]  {$g$};
				\draw [->,>=latex] (E) -- (A)   node[midway,below,rotate=0]  {$ $ $ $ $\gamma$};
				\draw [->,>=latex,dotted] (E) -- (B)   node[midway,below,rotate=0]  {$ $ $ $ $h$};
			\end{scope}
		\end{tikzpicture}
	\end{center}
	
	est commutatif).
\end{defn}
\begin{defn}\cite{Ko}.
	Un $A$-module $M$ est dit quasi-projectif (resp, quasi-injectif) si
	pour tout homomorphisme surjectif (resp, injectif) $g$ de M vers N
	(resp, de N vers M) et pour tout homomorphisme, $\gamma$ de M (resp,
	N) vers N (resp, vers M), il existe un endomorphisme $h$ de M tel
	que : $\gamma = gh$ (resp, $\gamma = hg$) (i.e., il existe $h:M
	\rightarrow M$ tel que le diagramme 
	
	\begin{center}
		\begin{tikzpicture}
			\begin{scope}[xscale=2.5,yscale=2.5]
				\node (A) at (1,0)   {$N$};
				\node (B) at (0,0){$M$};
				\node(C) at (2,0) {$0$};
				\node (E)   at (1,1) {$M$};
				\draw [->,>=latex] (A) -- (C)   node[midway,below,rotate=0]  {$ $};
				\draw [->,>=latex] (B) -- (A)   node[midway,above,rotate=0]  {$g$};
				\draw [->,>=latex] (E) -- (A)   node[midway,below,rotate=0]  {$ $ $ $ $\gamma$};
				\draw [->,>=latex,dotted] (E) -- (B)   node[midway,below,rotate=0]  {$ $ $ $ $h$};
			\end{scope}
		\end{tikzpicture}
	\end{center}
	
	resp,
	
	\begin{center}
		\begin{tikzpicture}
			\begin{scope}[xscale=2.5,yscale=2.5]
				\node (A) at (0,1)   {$N$};
				\node (B) at (1,1){$M$};
				\node(C) at (-1,1) {$0$};
				
				\node (E)   at (0,0) {$M$};
				\draw [->,>=latex] (C) -- (A)   node[midway,below,rotate=0]  {$ $};
				
				\draw [->,>=latex] (A) -- (B)   node[midway,above,rotate=0]  {$g$};
				\draw [->,>=latex] (A) -- (E)   node[midway,below,rotate=0]  {$ $ $ $ $\gamma$};

				\draw [->,>=latex,dotted] (B) -- (E)   node[midway,below,rotate=0]  {$ $ $ $ $h$};
				
			\end{scope}
		\end{tikzpicture}
	\end{center}
	
	est commutatif).
\end{defn}
\begin{prop}\cite{Sc}.
	Soit $M$ un $A$-module injectif, si $M$ est Hopfien alors il est
	co-Hopfien.
\end{prop}
\begin{prop}\cite{Sc}.
	Soit $M$ un $A$-module projectif, si $M$ est co-Hopfien alors il est
	Hopfien.
\end{prop}
\begin{prop}\cite{Mom}.
	Soit $M$ un $A$-module quasi-injectif, alors $M$ est Dedekind Fini
	si et seulement si $M$ est co-Hopfien.
\end{prop}
\begin{prop}\cite{Gh}.
	Soit $M$ un $A$-module quasi-projectif, alors $M$ est Dedekind Fini
	si et seulement si $M$ est Hopfien.
\end{prop}

\section{Quelques généralisations des modules superflus}
\begin{defn}\cite{Le}.
	Soit $M$ un $A$-module et soit $N$ un sous-module de $M$, on dit que $N$
	est superflu dans $M$ ($N\ll M$) si $N+L=M \Rightarrow L=M$, pour
	tout sous-module $L$ de $M$.
\end{defn}

Un module $M$ est dit creux si tout sous-module propre de $M$ est
superflu.

\begin{defn}\cite{Rotman}.
	L'enveloppe injective $E(M)$ d'un module $M$ est l'extension essentielle maximale de $M$.
\end{defn}
\begin{prop}\cite{Le}.
	Soit $M$ un $A$-module, alors $M$ est un module superflu si et
	seulement si $M$ est superflu dans son enveloppe injective $E(M)$.
\end{prop}

\begin{lem}\cite{Wi}
	Soient $M$, $N$ et $L$ des modules. Alors les deux épimorphismes
	$f : M \rightarrow N$ et $g : N \rightarrow L$ sont superflus si et
	seulement si $gf$ est superflu.
\end{lem}
\begin{defn}\cite{Go}.
	Soit $M$ un $A$-module, on appelle sous-module singulier de $M$,
	l'ensemble $Z(M)$ des éléments $x$ de $M$ tel que $Ann(x)$
	soit un idéal essentiel dans $A$. Un $A$-module $M$ est dit
	singulier (resp, non singulier) si $Z(M)=M$ (resp., $Z(M)=0)$.
	
\end{defn}

Si $N$ est un sous-module essentiel de $M$ alors $M/N$ est singulier,
mais la réciproque n'est pas en général vraie comme le montre
l'exemple suivant, soit $M=\mathbb{Z}/2\mathbb{Z}$ et $N=0$. $M/N$ est
singulier mais $N$ n'est pas essentiel dans $M$.

\begin{prop}\cite{Gha}.
	Soit $M$ un $A$-module injectif non singulier, alors $M$ est Hopfien
	si et seulement si $M$ est co-Hopfien.
\end{prop}
\begin{defn}\cite{Zh}.
	Soit $M$ un $A$-module et soit $N$ un sous-module de $M$, on dit que $N$
	est $\delta$-superflu dans $M$ ($N\ll_{\delta} M$) si $N+L=M$ tel
	que $M/L$ est singulier implique $L=M$, pour tout sous-module $L$ de
	$M$.
\end{defn}
\begin{lem}\label{lem:delta}\cite{Zh}. Soit $M$ un $A$-module.
	
	\begin{enumerate}
		
		\item Soient $K\leq B\leq M$. Alors $B\ll_{\delta}M$ si et seulement si $K\ll_{\delta}M$ et $B/K\ll_{\delta}M/K$.
		
		\item Soient $K$ et $B$ deux sous-modules de $M$, alors $K+B\ll_{\delta} M$ si et seulement si $K\ll_{\delta} M$ et $B\ll_{\delta} M$.
		
		\item Soient $K$ et $B$ deux sous-modules de $M$ avec $K\leq B$, si $K\ll_{\delta} B$,
		alors $K\ll_{\delta} M$.
		
		\item Soit $f :M \rightarrow N$ un homomorphisme tel que $K\ll_{\delta}
		M$, alors $f(K)\ll_{\delta} N$.
		
		\item Soit $M=M_{1}\oplus M_{2}$ un $A$-module et soient $A_{1}\leq M_{1}$
		et $A_{2}\leq M_{2}$, alors $A_{1}\oplus
		A_{2}\ll_{\delta}M_{1}\oplus M_{2}$ si et seulement si
		$A_{1}\ll_{\delta} M_{1}$ et $A_{2}\ll_{\delta} M_{2}$.
		
	\end{enumerate}
\end{lem}


\begin{defn}\cite{Oz}.
	Soit $M$ un $A$-module, on appelle sous-module cosingulier de $M$,
	l'ensemble $Z^{*}(M)$ des éléments $m$ de $M$ tel que $mA$
	soit un module superflu. Un $A$-module $M$ est dit cosingulier
	(resp, non cosingulier) si $Z^{*}(M)=M$ (resp., $Z^{*}(M)=0)$.
\end{defn}

\begin{lem}\label{lem:Ozc}\cite{Ozc}. Soit $M$ un $A$-module.
	
	\begin{enumerate}
		
		\item Si $M$ est superflu alors $Z^{*}(M)=M$.
		
		\item Si $M$ est semi simple injectif alors $Z^{*}(M)=0$

	\end{enumerate}
\end{lem}

\begin{lem}\label{lem:co}\cite{Was}.
	Soit $f :M\rightarrow N$ un homomorphisme et soit $A$ un sous-module
	de $M$ tel que $M/A$ est cosingulier, alors $f(M)/f(A)$ est
	cosingulier.
	
\end{lem}

\begin{defn}\cite{Was}.
	Soit $M$ un $A$-module et soit $N$ un sous-module de $M$, on dit que $N$
	est $\mu$-superflu dans $M$ ($N\ll_{\mu} M$) si $N+L=M$ tel que
	$M/L$ est cosingulier implique $L=M$, pour tout sous-module $L$ de $M$.
\end{defn}
\begin{lem}\label{lem:mu}\cite{Was}. Soit $M$ un $A$-module.
	
	\begin{enumerate}
		
		\item Soient $K\leq B\leq M$. Alors $B\ll_{\mu}M$ si et seulement si $K\ll_{\mu}M$ et $B/K\ll_{\mu}M/K$.
		
		\item Soient $K$ et $B$ deux sous-modules de $M$, alors $K+B\ll_{\mu} M$ si et seulement si $K\ll_{\mu} M$ et $B\ll_{\mu} M$.
		
		\item Soient $K$ et $B$ deux sous-modules de $M$ avec $K\leq B$, si $K\ll_{\mu} B$,
		alors $K\ll_{\mu} M$.
		
		\item Soit $f :M \rightarrow N$ un homomorphisme tel que $K\ll_{\mu}
		M$, alors $f(K)\ll_{\mu} N$.
		
		\item Soit $M=M_{1}\oplus M_{2}$ un $A$-module et soient $A_{1}\leq M_{1}$
		et $A_{2}\leq M_{2}$, alors $A_{1}\oplus A_{2}\ll_{\mu}M_{1}\oplus
		M_{2}$ si et seulement si $A_{1}\ll_{\mu} M_{1}$ et $A_{2}\ll_{\mu}
		M_{2}$.
		
	\end{enumerate}
\end{lem}

\begin{lem}\cite{Was}.\label{lem:mossa}
	Soient $M$ un $A$-module et $K\leq N$ deux sous-modules
	de $M$ , si $N$ est un facteur direct de $M$ et
	$K\ll_{\mu} M$, alors $K\ll_{\mu} N$.
\end{lem}
\begin{defn}\cite{Me}.
	Soit $M$ un $A$-module et soit $N$ un sous-module de $M$, on dit que $N$
	est $\gamma$-superflu dans $M$ ($N\ll_{\gamma} M$) si $N+L=M$ tel que
	$M/L$ est non cosingulier implique $L=M$, pour tout sous-module $L$ de $M$.
\end{defn}
\begin{lem}\label{lem:gamma}\cite{Me}. Soit $M$ un $A$-module.
	
	\begin{enumerate}
		
		\item Soient $K\leq B\leq M$. Alors $B\ll_{\gamma}M$ si et seulement si $K\ll_{\gamma}M$ et $B/K\ll_{\gamma}M/K$.

		\item Soient $K$ et $B$ deux sous-modules de $M$ avec $K\leq B$, si $K\ll_{\gamma} B$,
		alors $K\ll_{\gamma} M$.
		
		\item Soit $f :M \rightarrow N$ un épimorphisme tel que $K\ll_{\gamma}
		M$, alors $f(K)\ll_{\gamma} N$.
		
		\item Soit $M=M_{1}\oplus M_{2}$ un $A$-module et soient $A_{1}\leq M_{1}$
		et $A_{2}\leq M_{2}$, alors $A_{1}\oplus A_{2}\ll_{\gamma}M_{1}\oplus
		M_{2}$ si et seulement si $A_{1}\ll_{\gamma} M_{1}$ et $A_{2}\ll_{\gamma}
		M_{2}$.
		
	\end{enumerate}
\end{lem}

\begin{defn}\cite{An}.
	Soit $M$ un $A$-module.
	
	$\mathrm{Rad}(M) = \bigcap \{ K\leq M / K $ est maximal dans $M\}$=
	$\sum \{K\leq M / K$ est superflu dans $M\}$
	
\end{defn}
\begin{defn}\cite{Ka}.
	Soit $M$ un $A$-module et soit $N$ un sous-module de $M$, on dit que $N$
	est Jacobson-superflu dans $M$ ($N\ll_{J} M$) si $N+L=M$ tel que
	$\mathrm{Rad}(M/L)=M/L$ implique $L=M$, pour tout sous-module $L$ de $M$.
\end{defn}
\begin{lem}\label{lem:J}\cite{Ka}. Soit $M$ un $A$-module.
	
	\begin{enumerate}
		
		\item Soient $K\leq B\leq M$. Alors $B\ll_{J}M$ si et seulement si $K\ll_{J}M$ et $B/K\ll_{J}M/K$.
		
		\item Soient $K$ et $B$ deux sous-modules de $M$, alors $K+B\ll_{J} M$ si et seulement si $K\ll_{J} M$ et $B\ll_{J} M$.
		
		\item Soient $K$ et $B$ deux sous-modules de $M$ avec $K\leq B$, si $K\ll_{J} B$,
		alors $K\ll_{J} M$.
		
		\item Soit $f :M \rightarrow N$ un homomorphisme tel que $K\ll_{J}
		M$, alors $f(K)\ll_{J} N$.
		
		\item Soit $M=M_{1}\oplus M_{2}$ un $A$-module et soient $A_{1}\leq M_{1}$
		et $A_{2}\leq M_{2}$, alors $A_{1}\oplus A_{2}\ll_{J}M_{1}\oplus
		M_{2}$ si et seulement si $A_{1}\ll_{J} M_{1}$ et $A_{2}\ll_{J}
		M_{2}$.
		
	\end{enumerate}
\end{lem}


\section{Modules Hopfien généralisé et faiblement co-Hopfien}

\begin{defn}\cite{Gh}.
	Un $A$-module $M$ est appelé Hopfien généralisé, si
	pour tout endomorphisme surjectif $f$ de $M$, le noyau de $f$ est
	superflu dans $M$ $(Ker f \ll M)$.
\end{defn}
\begin{cor}\cite{Gh}.
	Soit $M$ un $A$-module quasi-projectif, alors les assertions
	suivantes sont équivalentes:
	\begin{enumerate}
		
		\item $M$ est Hopfien.
		
		\item $M$ est Hopfien généralisé.
		
		\item $M$ est Dedekind Fini.
		
	\end{enumerate}
	
\end{cor}
\begin{defn}\cite{Wi}.
	Soient $M$ un $A$-module et $N$ un sous-module de $M$, on dit que $N$
	est essentiel de $M$ ($N\leq^{e} M$) si $N\cap L=0 \Rightarrow L=0$,
	pour tout sous-module $L$ de $M$.
\end{defn}

Un module $M$ est dit uniforme si tout sous-module non nul de $M$ est
essentiel.
\begin{defn}\cite{Ha}.
	Un $A$-module $M$ est appelé faiblement co-Hopfien, si pour tout
	endomorphisme injectif $f$ de $M$, l'image de $f$ est essentiel de
	$M$ $(Im f \leq^{e} M)$.
\end{defn}
\begin{cor}\cite{Ha}.
	Soit $M$ un $A$-module quasi-injectif, alors les assertions
	suivantes sont équivalentes:
	\begin{enumerate}
		
		\item $M$ est co-Hopfien.
		
		\item $M$ est faiblement co-Hopfien.
		
		\item $M$ est Dedekind Fini
		
	\end{enumerate}
	
\end{cor}
\begin{rems}\cite{Gh}\cite{Ha}. Soit $M$ un $A$-module.
	\begin{enumerate}

		\item Si $M$ vérifie la CCD sur les sous-modules non essentiels,
		alors $M$ est faiblement co-Hopfien.
		
		\item Si $M$ vérifie la CCA sur les sous-modules non superflus,
		alors $M$ est Hopfien généralisé.
		
	\end{enumerate}
\end{rems}

\section{Modules faiblement Hopfien et co-Hopfien généralisé}
\begin{defn}\cite{Wa}.
	Un $A$-module $M$ est dit faiblement Hopfien si tout endomorphisme
	surjectif superflu $f$ $(Ker f \ll M)$ de $M$ est bijectif.
\end{defn}
\begin{defn}\cite{Wa}.
	Un $A$-module $M$ est dit co-Hopfien généralisé si tout
	endomorphisme injectif essentiel $f$ $(Im f \leq^{e} M)$ de $M$ est
	bijectif.
\end{defn}
\begin{rems}\cite{Wa}. Soit $M$ un $A$-module.
	\begin{enumerate}
		
		\item $M$ est co-Hopfien si et seulement si $M$ est faiblement co-Hopfien
		et co-Hopfien généralisé.
		
		\item $M$ est Hopfien si et seulement si $M$ est faiblement Hopfien
		et Hopfien généralisé.
		
		\item Si $M$ vérifie la CCD sur les sous-modules essentiels
		alors $M$ est co-Hopfien généralisé.
		
		\item Si $M$ vérifie la CCA sur les sous-modules superflus
		alors $M$ est faiblement Hopfien.
		
	\end{enumerate}
\end{rems}



\begin{thm}\cite{Pal}
	Soit $M$ un module quasi-projectif et soit $N$ un sous-module totalement
	invariant superflu dans $M$, si $M$ est faiblement Hopfien alors
	$M/N$ est faiblement Hopfien.
\end{thm}

\begin{cor}\cite{Pal}
	Soit $M$ un module quasi-projectif de type fini, si $M$ est
	faiblement Hopfien alors $M/Rad(M)$ est faiblement Hopfien.
\end{cor}

\begin{prop}\cite{Pal}
	Soit $M$ un module quasi-projectif, si $M$ est co-Hopfien alors il
	est faiblement Hopfien.
\end{prop}

\begin{prop}\cite{Pal}
	Soit $M$ un module quasi-injectif, si $M$ est Hopfien alors il est
	co-Hopfien généralisé.
\end{prop}

\section{Modules fortement Hopfien et fortement co-Hopfien}
\begin{defn}\cite{Hm}.
	Un $A$-module $M$ est appelé fortement Hopfien, si pour tout
	endomorphisme $f$ de $M$ la suite croissante: $Ker f \subseteq Ker
	f^{2} \subseteq  ... \subseteq Ker f^{n} \subseteq ...$ est
	stationnaire.
	
\end{defn}
\begin{prop}\cite{Hm}.
	Soit $M$ un $A$-module, alors les assertions suivantes sont
	équivalentes:
	\begin{enumerate}
		
		\item $M$ est fortement Hopfien.
		
		\item Pour tout endomorphisme $f$ de $M$, il existe $n\geq 1$ tel que $Ker f^{n}=Ker f^{n+1}$.
		
		\item Pour tout endomorphisme $f$ de $M$, il existe $n\geq 1$ tel que $Ker f^{n}\cap Im f^{n}=(0)$.
		
	\end{enumerate}
	
\end{prop}
\begin{defn}\cite{Hm}.
	Un $A$-module $M$ est appelé fortement co-Hopfien, si pour tout
	endomorphisme $f$ de $M$ la suite décroissante: $Im f \supseteq
	Im f^{2} \supseteq  ... \supseteq Im f^{n} \supseteq ...$ est
	stationnaire.
	
\end{defn}
\begin{prop}\cite{Hm}.
	Soit $M$ un $A$-module, alors les assertions suivantes sont
	équivalentes:
	\begin{enumerate}
		
		\item $M$ est fortement co-Hopfien.
		
		\item Pour tout endomorphisme $f$ de $M$, il existe $n\geq 1$ tel que $Im f^{n}=Im f^{n+1}$.
		
		\item Pour tout endomorphisme $f$ de $M$, il existe $n\geq 1$ tel que $M= Ker f^{n}+ Im f^{n}$.
		
	\end{enumerate}
	
\end{prop}
\begin{defn}\cite{Ar}.
	Un $A$-module $M$ est dit un module de Fitting si pour tout
	endomorphisme $f$ de $M$, il existe un entier $n\geq 1$ tel que
	$M=Ker f^{n}\oplus Im f^{n}$.
	
\end{defn}
\begin{rems}\cite{Hm}.
	Soit $M$ un $A$-module:
	\begin{enumerate}
		
		\item Tout module noethérien (resp, artinien) est fortement Hopfien (resp,
		fortement co-Hopfien).
		
		\item Tout module fortement Hopfien (resp, fortement co-Hopfien) est Hopfien (resp,
		co-Hopfien).
		
		\item Tout module fortement Hopfien (resp, fortement co-Hopfien) est
		Dedekind Fini.
		
		\item $M$ est fortement Hopfien et fortement co-Hopfien si et
		seulement si $M$ est un module de Fitting.
		
		\item Tout module de longueur finie est de Fitting.
		
		\item Tout module de Fitting est Hopfien et co-Hopfien.
	\end{enumerate}
	
\end{rems}
\begin{ex}\cite{Hm}\cite{Yl}.
	Il existe un $\mathbb{Z}$-module Hopfien et co-Hopfien et qui n'est
	ni fortement Hopfien ni fortement co-Hopfien. En effet: Soit
	$(p_{n})_{n\geq 1}$ une suite de nombres premiers tels que
	$p_{1}<p_{2}< ... <p_{n}<...$. Posons
	$M_{n}=\mathbb{Z}/p_{n}^{n}\mathbb{Z}$ et considérons
	$M=\oplus_{n\geq 1}M_{n}$. $M$ est co-Hopfien (resp, Hopfien), par
	contre $M$ n'est ni fortement Hopfien ni fortement co-Hopfien.
	D'autre part, soit $P$ l'ensemble des nombres premiers. Le
	$\mathbb{Z}$-module $M=\oplus_{p\in P}\mathbb{Z}_{p}$ est fortement
	Hopfien et fortement co-Hopfien mais n'est ni artinien ni
	noethérien.
	
\end{ex}

\begin{thm}\cite{Hm}.
	Soit $M$ un $A$-module:
	\begin{enumerate}
		
		\item Si $M$ est quasi-projectif et fortement co-Hopfien alors $M$
		est fortement Hopfien.
		
		\item Si $M$ est quasi-injectif et fortement Hopfien alors $M$
		est fortement co-Hopfien.
		
	\end{enumerate}
	
\end{thm}

\section{Modules semi Hopfien et semi co-Hopfien}
\begin{defn}\cite{Pi}.
	Un $A$-module $M$ est appelé semi Hopfien, si pour tout
	endomorphisme surjectif $f$ de $M$, le noyau de $f$ est un facteur
	direct de $M$.
	
\end{defn}

\begin{exs}\cite{Pi}.
	\begin{enumerate}
		
		\item Tout module semi simple est semi Hopfien.
		
		\item D'après \cite[Théorème 16(ii)]{Hi}, un espace vectoriel V sur un corps F est
		Hopfien si et seulement s'il est de dimension finie. Donc un espace vectoriel de dimension infinie
		est semi Hopfien, mais il n'est pas Hopfien.
		
		\item Tout module vérifiant D2 est semi Hopfien. (On dit qu'un module $M$ vérifie D2 si tout sous-module $N$ tel que $M/N$ est isomorphe à un facteur direct
		de $M$ est un facteur direct de $M$).
		
		\item Tout module quasi-projectif est semi Hopfien.

	\end{enumerate}
	
\end{exs}

\begin{prop}\cite{Pi}.
	Soit $M$ un $A$-module semi-Hopfien, si $M$ est Dedekind Fini alors
	il est Hopfien.
\end{prop}
\begin{defn}\cite{Pi}.
	Un $A$-module $M$ est appelé semi co-Hopfien, si pour tout
	endomorphisme injectif $f$ de $M$, l'image de $f$ est un facteur
	direct de $M$.
	
\end{defn}

\begin{prop}\cite{Pi}.
	Soit $M$ un $A$-module semi-co-Hopfien, si $M$ est Dedekind Fini
	alors il est co-Hopfien.
\end{prop}

\begin{rem}
	Les anneaux fortement $\pi$-réguliers à gauche et à
	droite ont été introduits par Kaplansky \cite{Kap}, Azumaya
	a montré en 1954 qu'un anneau $A$ est fortement
	$\pi$-régulier si pour tout $a\in A$ il existe $m\in \mathbb{N}$
	et $c\in A$ tel que $ac=ca$ et $a^{m}=ca^{m+1}$ \cite{Az}.
	Dischinger a montré en 1976 que la propriété fortement
	$\pi$-régularité est symétrique \cite{Di}.
\end{rem}

\begin{ex}
	
	D'après \cite[Remarque 2.16(3)]{Hm}, l'anneau $A=\prod _{n\geq1}
	\mathbb{Z}/2^{n}\mathbb{Z}$ est Hopfien (tout anneau commutatif est
	Hopfien) mais n'est pas fortement Hopfien. Puisque tout anneau
	Hopfien est semi Hopfien, alors l'anneau $\prod _{n\geq1}
	\mathbb{Z}/2^{n}\mathbb{Z}$ est semi Hopfien mais n'est pas
	fortement Hopfien.
\end{ex}
\begin{thm}\cite{Jor}.
	Soit $M$ un $A$-module, alors:
	\begin{enumerate}
		
		\item Si $M$ est semi Hopfien et fortement co-Hopfien, alors
		$End_{A}(M)$ est fortement $\pi$-régulier.
		
		\item Si M est semi co-Hopfien et fortement Hopfien, alors
		$End_{A}(M)$ est fortement $\pi$-régulier.
	\end{enumerate}
	
\end{thm}

\begin{cor}\cite{Jor}
	Tout module semi Hopfien et fortement co-Hopfien ou semi co-Hopfien et fortement
	Hopfien est un module de Fitting.
\end{cor}
Le résultat suivant présente un analogue du
théorème de Hopkins-Levitzki.
\begin{cor}\cite{Jor}
	Soit $M$ un $A$-module, alors:
	
	\begin{enumerate}
		
		\item Si $M$ est semi Hopfien et fortement co-Hopfien, alors $M$ est
		fortement Hopfien.
		
		\item Si M est semi co-Hopfien et fortement Hopfien, alors $M$ est
		fortement co-Hopfien.
		
	\end{enumerate}
\end{cor}


Il est facile de voir que tout module Hopfien est semi Hopfien, mais
la réciproque n'est pas vraie en général comme le montre l'exemple suivant:

\begin{ex}\cite{Jor}
	D'après \cite[Théorème 16(ii)]{Hi}, un espace vectoriel
	$V$ sur un corps $F$ est Hopfien si et seulement s'il est de dimension
	finie. Alors un espace vectoriel de dimension infinie sur un corps
	est semi Hopfien mais n'est pas Hopfien.
	
\end{ex}
\begin{prop}\cite{Jor}
	Soit $M$ un module semi Hopfien, si $M$ est indécomposable alors il est
	Hopfien.
\end{prop}


\begin{thm}\cite{Jor}
	Soit $M$ un $A$-module, alors:
	
	(1) Si $M$ est semi Hopfien et co-Hopfien, alors $M$ est Hopfien.
	
	(2) Si $M$ est semi co-Hopfien et Hopfien, alors $M$ est co-Hopfien.

\end{thm}


\begin{defn}\cite{Ca}.
	Un $A$-module $M$ est dit quasi-principalement projectif si pour
	tout endomorphisme $f$ de $M$ et pour tout homomorphisme $g$ de M
	vers $f(M)$, il existe un endomorphisme $h$ de $M$ tel que : $g = fh$
	
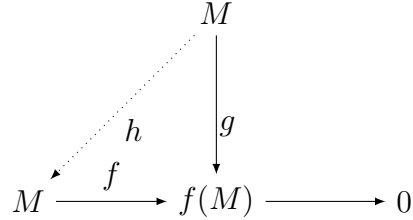
\begin{center}
		\begin{tikzpicture}
			\begin{scope}[xscale=2.5,yscale=2.5]
				\node (A) at (1,0)   {$f(M)$};
				\node (B) at (0,0){$M$};
				\node(C) at (2,0) {$0$};
				
				\node (E)   at (1,1) {$M$};
				\draw [->,>=latex] (A) -- (C)   node[midway,below,rotate=0]  {$ $};
				
				\draw [->,>=latex] (B) -- (A)   node[midway,above,rotate=0]  {$f$};
				\draw [->,>=latex] (E) -- (A)   node[midway,below,rotate=0]  {$ $ $ $ $g$};
				\draw [->,>=latex,dotted] (E) -- (B)   node[midway,below,rotate=0]  {$ $ $ $ $h$};
			\end{scope}
		\end{tikzpicture}
	\end{center}
\end{defn}

Puisque tout $A$-module quasi-principalement projectif est semi
Hopfien d'après \cite[Proposition 3.2]{Kum}, alors il est facile
de voir le corollaire suivant.

\begin{cor}\cite{Jor}
	Soit $M$ un module quasi-principalement projectif, si $M$ est co-Hopfien
	alors il est Hopfien.
\end{cor}


\begin{defn}\cite{Ni}.
	Un $A$-module $M$ est dit quasi-principalement injectif si pour tout
	endomorphisme non nul $f$ de M et pour tout homomorphisme $g$ de
	$f(M)$ vers $M$, il existe un endomorphisme $h$ de $M$ tel que : $g =
	hf$
	
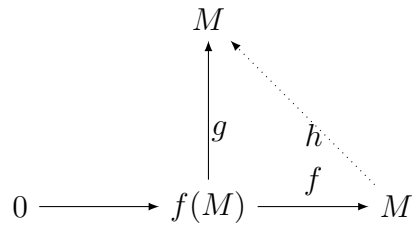
\begin{center}
		\begin{tikzpicture}
			\begin{scope}[xscale=2.5,yscale=2.5]
				\node (A) at (1,0)   {$f(M)$};
				\node (B) at (2,0){$M$};
				\node(C) at (0,0) {$0$};
				
				\node (E)   at (1,1) {$M$};
				\draw [->,>=latex] (C) -- (A)   node[midway,below,rotate=0]  {$ $};
				
				\draw [->,>=latex] (A) -- (B)   node[midway,above,rotate=0]  {$f$};
				\draw [->,>=latex] (A) -- (E)   node[midway,below,rotate=0]  {$ $ $ $ $g$};
				\draw [->,>=latex,dotted] (B) -- (E)   node[midway,below,rotate=0]  {$ $ $ $ $h$};
			\end{scope}
		\end{tikzpicture}
	\end{center}
\end{defn}

Puisque tout $A$-module quasi-principalement injectif est semi
co-Hopfien d'après \cite[Proposition 3.1]{Kum}, alors il est
facile de voir le corollaire suivant.

\begin{cor}\cite{Jor}
	
	Soit $M$ un module quasi-principalement injectif, si $M$ est Hopfien
	alors il est co-Hopfien.
	
\end{cor}

Il est facile de voir que tout module Hopfien est Hopfien
généralisé, mais la réciproque n'est pas toujours vraie comme le montre l'exemple suivant.

\begin{ex}\cite[exemple 1.7]{Gh}.
	Soit $G=\mathbb{Z}_{p^{\infty}}$. Puisque dans $G$ tout sous-groupe
	propre est superflu, donc $G$ est un groupe abélien Hopfien
	généralisé. Mais $G$ n'est pas Hopfien puisque la
	multiplication par $p$ induit un épimorphisme de $G$ qui n'est
	pas un isomorphisme.
	
\end{ex}
\begin{prop}\cite{Jor}
	Soit $M$ un module semi Hopfien. Alors les assertions suivantes sont
	équivalentes:
	
	(1) $M$ est Hopfien.
	
	(2) $M$ est Hopfien généralisé.
	
\end{prop}


\begin{prop}\cite{Jor}
	Soit $M$ un module semi co-Hopfien. Alors les assertions suivantes
	sont équivalentes:
	
	(1) $M$ est co-Hopfien.
	
	(2) $M$ est faiblement co-Hopfien
	
\end{prop}

\section{Modules $\mu$-Hopfiens}
\begin{defn}\cite{Fer}
	
	Un $A$-module $M$ est appelé $\mu$-Hopfien, si pour tout
	endomorphisme surjectif $f$ de $M$, le noyau de $f$ est
	$\mu$-superflu dans $M$ $(Ker f \ll_{\mu} M)$.
	
\end{defn}
\begin{lem}\label{lem:abd}\cite{Fer}
	Soit $M$ un $A$-module et soit $N$ un sous-module de $M$, alors les
	assertions suivantes sont équivalentes:
	\begin{enumerate}
		
		\item $N\ll_{\mu} M$.
		
		\item Si $X+N=M$, alors $X$ est un facteur direct de $M$ avec $M/X$ est un module semi simple
		injectif.
		
	\end{enumerate}
	
\end{lem}

Le résultat suivant présente une caractérisation des
modules $\mu$ -Hopfiens.

\begin{thm}\label{thm:ps}\cite{Fer} Soit $M$ un $A$-module, alors les
	assertions suivantes sont équivalentes:
	\begin{enumerate}
		
		\item $M$ est $\mu$-Hopfien.
		
		\item Pour tout endomorphisme surjectif $f$ de $M$, si $N\ll_{\mu} M$,
		alors $f^{-1}(N) \ll_{\mu}M$.
		
		\item Si $N\leq M$ et s'il existe un épimorphisme $M/N\rightarrow M$,
		alors $N\ll_{\mu} M$.
		
		\item Si $M/N$ est non nul et cosingulier pour tout $N\leq M$ et si $f$
		est un endomorphisme surjectif de $M$ alors $f(N)\neq M$.

		\item Il existe un sous-module $\mu$-superflu totalement invariant $N$ de $M$ tel que $M/N$ est $\mu$-Hopfien.
		
		\item Pour tout module $X$ tel qu'il existe un épimorphisme $M\rightarrow M\oplus X$, alors $X$ est semi simple injectif.

	\end{enumerate}
\end{thm}


L'exemple suivant montre que la classe des modules Hopfiens est une
sous-classe propre des modules $\mu$-Hopfiens.

\begin{ex}\label{ex:4s}\cite{Fer}
	
	Soit $G=\mathbb{Z}_{p^{\infty}}$. Puisque dans $G$ chaque
	sous-groupe propre est $\mu$-superflu. Donc $G$ est un groupe
	$\mu$-Hopfien. Mais $G$ n'est pas Hopfien car la multiplication par
	$p$ induit un épimorphisme de $G$ qui n'est pas un isomorphisme.
\end{ex}



\begin{thm}\cite{Fer}\label{thm:pp} Soit $M$ un module (quasi-)projectif et soit $f\in \mathrm{End}(M)$.
	Alors les assertions suivantes sont équivalentes:
	\begin{enumerate}
		
		\item $M$ est $\mu$-Hopfien.
		
		\item Si $f$ est un épimorphisme, alors $Ker(f)$ est semi simple injectif.

	\end{enumerate}
\end{thm}



\begin{thm}\cite{Fer}\label{thm:ppp} Soit $A$ un anneau. Alors les assertions suivantes sont équivalentes:
	\begin{enumerate}
		
		\item Tout $A$-module est $\mu$-Hopfien.
		
		\item Tout $A$-module projectif est $\mu$-Hopfien.
		
		\item Tout $A$-module libre est $\mu$-Hopfien.
		
		\item $A$ est semi simple.
		
	\end{enumerate}
\end{thm}


Il est clair que tout module Hopfien généralisé est
$\mu$-Hopfien. L'exemple suivant montre que la réciproque n'est pas vraie
en général. De plus, cet exemple montre aussi qu'un module
$\mu$-Hopfien n'est pas nécessairement Dedekind Fini.

\begin{ex}\label{ex:exs}\cite{Fer}
	
	Soit $A$ un anneau semi simple. Alors d'après le
	théorème~\ref{thm:ppp}, $M=A^{(\mathbb{N})}$ est un $A$-module
	$\mu$-Hopfien. Puisque $A^{(\mathbb{N})}\cong
	A^{(\mathbb{N})}\oplus A^{(\mathbb{N})}$ et $A^{(\mathbb{N})}\neq
	0$, $M$ n'est pas un module Hopfien généralisé (Dedekind
	Fini) (voir \cite[Corollaire 1.4]{Gh}).
\end{ex}


\begin{prop}\cite{Fer}
	Soit $N$ un sous-module totalement invariant de $M$ tel que $M/N$
	est Hopfien. Si $N$ est $\mu$-Hopfien alors $M$ l'est.
\end{prop}


\begin{prop}\label{pr:4}\cite{Fer}
	Soit $M$ un $A$-module. Si $M$ vérifie CCA sur les sous modules
	non $\mu$-superflus alors il est $\mu$-Hopfien.
\end{prop}


\begin{prop}\cite{Fer}
	Soit $M$ un $A$-module vérifie la propriété suivante, pour
	tout endomorphisme $f$ de $M$ il existe un entier $n\geq 1$ tel que
	$Kerf^{n}\cap Imf^{n}\ll_{\mu}M$. Alors $M$ est $\mu$-Hopfien.
	
\end{prop}


\begin{thm}\label{thm:mor}\cite{Fer}
	La propriété $\mu$-Hopfien est préservée par l'équivalence de Morita.
\end{thm}



\begin{prop}\label{prop:100}\cite{Fer}
	Tout facteur direct d'un module $\mu$-Hopfien $M$ est $\mu$-Hopfien.
\end{prop}


\begin{prop}\label{prop:11}\cite{Fer}
	Soit $M = M_{1}\oplus M_{2}$ un $A$-module. Si pour tout $i\in
	\{1,2\}$, $M_{i}$ est un sous-module totalement invariant de $M$,
	alors $M$ est $\mu$-Hopfien si et seulement si $M_{i}$ est
	$\mu$-Hopfien pour tout $i\in \{1,2\}$.
\end{prop}


\begin{defn}\cite{Fer}
	Soient $M$ et $N$ deux $A$-modules. $M$ est appelé $\mu$-Hopfien
	relatif à $N$, si pour tout épimorphisme $f : M\rightarrow
	N$, $Ker(f)\ll_{\mu}M$.
\end{defn}

\section{Modules $\delta$- faiblement Hopfiens}

\begin{defn}\cite{Kho}
	
	Un $A$-module $M$ est dit $\delta$-faiblement Hopfien si tout
	endomorphisme surjectif $\delta$-superflu $(Ker f \ll_{\delta} M)$
	de $M$ est bijectif.
	
\end{defn}

\begin{ex}\label{ex:z}\cite{Kho}
	Il existe un épimorphisme $\delta$-superflu qui n'est pas un
	isomorphisme. Soit $G = \mathbb{Z}_{p^{\infty}}$, comme dans $G$ tout sous
	groupe propre est $\delta$-superflu (car dans $G$ tout sous groupe
	propre est superflu), donc tout endomorphisme surjectif de $G$ est
	$\delta$-superflu, mais la multiplication par $p$ induit un
	épimorphisme de $G$ qui n'est pas un isomorphisme.
\end{ex}


\begin{lem}\label{lem:3}\cite{Kho}
	Soit $M$ un $A$-module. Alors les assertions suivantes sont
	équivalentes:
	
	\begin{enumerate}
		
		\item $M$ est $\delta$-faiblement Hopfien.
		
		\item Pour tout sous-module $\delta$-superflu $K$ de $M$, $M/K \cong M$ si et seulement si $K = 0$.
		
	\end{enumerate}
\end{lem}


\begin{prop}\label{prop:4ab}\cite{Kho}
	Soit $M$ un module $\delta$-faiblement Hopfien. Si $M\cong M\oplus
	N$ pour certain module projectif semi simple $N$, alors $N=0$. De
	plus, si $M$ est projectif, alors la réciproque est vraie.
\end{prop}

\begin{prop}\label{prop:4aa}\cite{Kho}
	Soit $A$ un anneau semi simple artinien. Alors un $A$-module libre
	$F$ est $\delta$-faiblement Hopfien si et seulement s'il est de rang
	fini.
\end{prop}

Le résultat suivant présente une caractérisation des modules
projectifs $\delta$-faiblement Hopfiens.

\begin{thm}\label{thm:6a}\cite{Kho}
	Soit $M$ un module projectif et $f \in End(M)$, alors les assertions
	suivantes sont équivalentes:
	\begin{enumerate}
		
		\item $M$ est $\delta$-faiblement Hopfien.
		
		\item Si $f$ est inversible à droite et $Ker(f)$ est semi simple, alors $f$ est inversible à gauche.
		
		\item Si $f$ est inversible à droite et $Ker(f)\ll_{\delta}M$, alors $f$ est inversible à gauche.
		
		\item Si $f$ admet un inverse à droite $g$ et $(1-gf)M\ll_{\delta}M$, alors $f$ est inversible à gauche.
		
		\item Si $f$ est surjectif et $Ker(f)$ est semi simple projectif, alors $f$ est inversible à gauche.

	\end{enumerate}
\end{thm}

Le résultat suivant présente une caractérisation des anneaux dans
lesquels tout module quasi-projectif (projectif, libre) est
$\delta$-faiblement Hopfien.

\begin{thm}\label{thm:6aa}\cite{Kho}
	Soit $A$ un anneau, alors les assertions suivantes sont
	équivalentes:
	\begin{enumerate}
		
		\item Tout $A$-module quasi-projectif est $\delta$-faiblement Hopfien.
		
		\item Tout $A$-module projectif est $\delta$-faiblement Hopfien.
		
		\item Tout $A$-module libre est $\delta$-faiblement Hopfien.
		
		\item Tout idéal à droite maximal de $A$ est essentiel dans
		$A_{A}$.
		
		\item $A$ n'a pas de $A$-module non nul semi simple projectif.
		
		\item $\delta(A)=J(A)$.
		
	\end{enumerate}
\end{thm}

On dit qu'un anneau $A$ est $GV$-anneau à droite \cite{Ra}, si tout $A$-module simple est soit projectif, soit injectif. Il est clair qu'un anneau $A$ est $GV$-anneau
si et seulement si tout $A$-module simple singulier est injectif. Noter aussi que d'après \cite[Corollaire 3.3]{Ta}, un anneau $A$ est $GV$-anneau
si et seulement si et seulement si tout $A$-module superflu est projectif.

\begin{cor}\cite{Kho}
	Soit $A$ un $GV$-anneau. Alors tout $A$-module superflu indécomposable est $\delta$-faiblement Hopfien.		
\end{cor}


Il est clair que tout module $\delta$-faiblement Hopfien est
faiblement Hopfien. L'exemple suivant montre que la réciproque
n'est pas toujours vraie.

\begin{ex}\label{ex:2}\cite{Kho}
	D'après \cite[Example 3.13]{Wa} tout espace vectoriel de
	dimension infinie est faiblement Hopfien. Mais d'après \cite[Lemme 2.9]{Tr}, $M \ll_{\delta} M$ si $M$ est un
	$A$-module projectif semi simple, donc tout endomorphisme surjectif de $M$
	est $\delta$-superflu. Alors $M$ n'est pas $\delta$-faiblement
	Hopfien lorsque $M$ n'est pas Dedekind Fini, c'est le cas des
	espaces vectoriels de dimension infinie.
\end{ex}


\begin{ex}\label{ex:3}\cite{Kho}
	Soit $P$ l'ensemble de tous les nombres premiers et
	$\mathbb{Q}/\mathbb{Z}=\bigoplus_{p\in P}\mathbb{Z}_{p^{\infty}}$. Si $\bigoplus_{p\in P}\mathbb{Z}_{p^{\infty}}$ est un $\mathbb{Z}$-module $\delta$-faiblement Hopfien, alors $\mathbb{Z}_{p^{\infty}}$ est
	$\delta$-faiblement Hopfien d'après la proposition~\ref{prop:10}, contradiction avec l'exemple~\ref{ex:z}. Donc $\mathbb{Q}/\mathbb{Z}$ n'est pas
	$\delta$-faiblement Hopfien, mais $\mathbb{Q}$ est un $\mathbb{Z}$-module $\delta$-faiblement Hopfien.
	
\end{ex}

\begin{thm}\label{thm:6}\cite{Kho}
	Soit $M$ un module quasi-projectif uniforme, si $N$ est un
	sous-module non nul totalement invariant $\delta$-superflu de $M$ alors
	$M/N$ est Hopfien.
\end{thm}


\begin{defn}\cite{Zh}
	Soit $\mathfrak{S}$ la classe de tous les modules simples singuliers. Pour un
	module $M$, soit $\delta(M)=Rej_{M}(\mathfrak{S})=\cap \{N \subseteq M ; M/N \in
	\mathfrak{S}\}$.
\end{defn}


\begin{cor}	\cite{Kho}
	Soit $M$ un module quasi-projectif uniforme tel que $\delta(M)$ est $\delta$-superflu dans $M$. Alors $M/\delta(M)$ est $\delta$-faiblement Hopfien.
\end{cor}


\begin{prop}\label{prop:8}\cite{Kho}
	Soit $M$ un module quasi-projectif, si $M$ est co-Hopfien alors il
	est $\delta$-faiblement Hopfien.
\end{prop}



\begin{defn}\cite{Oza}
	Soit $M$ un $A$-module. On dit que $M$ est duo module si tout sous module de $M$ est totalement invariant.
\end{defn}

\begin{cor}\cite{Kho}
	Soit $M = M_{1}\oplus M_{2}$ un duo module. Alors $M$ est $\delta$-faiblement
	Hopfien si et seulement si $M_{1}$ et $M_{2}$ sont $\delta$-faiblement Hopfiens.
\end{cor}


Il est clair que tout module Hopfien est $\delta$-faiblement
Hopfien. L'exemple suivant montre que la réciproque n'est pas
toujours vraie, Ainsi il montre qu'un module $\delta$-faiblement
Hopfien n'est pas toujours Dedekind Fini.

\begin{ex}\label{ex:1}\cite{Kho}
	Si $M$ est un $A$-module singulier semi simple, alors le seul
	sous-module $\delta$-superflu de $M$ est zéro. Donc tout
	endomorphisme surjectif $\delta$-superflu de $M$ est injectif. Mais
	$M$ n'est pas Dedekind Fini et donc n'est pas Hopfien.
	
\end{ex}

\begin{thm}\label{thm:12}\cite{Kho}
	Soit $M$ un $A$-module vérifiant la CCA sur les sous-modules
	$\delta$-superflus. Alors $M$ est $\delta$-faiblement Hopfien.
\end{thm}


\section{Modules $\gamma$-Hopfiens}

\begin{defn}\cite{Alg}
	
	Un $A$-module $M$ est appelé $\gamma$-Hopfien, si pour tout
	endomorphisme surjectif $f$ de $M$, le noyau de $f$ est
	$\gamma$-superflu dans $M$ $(Ker f \ll_{\gamma} M)$.
	
\end{defn}
Le résultat suivant présente une caractérisation des
Modules $\gamma$ -Hopfiens.

\begin{thm}\label{thm:p}\cite{Alg} Soit $M$ un $A$-module, alors les
	assertions suivantes sont équivalentes:
	\begin{enumerate}
		
		\item $M$ est $\gamma$-Hopfien,
		
		\item Pour tout endomorphisme surjectif $f$ de $M$, si $N\ll_{\gamma} M$,
		alors $f^{-1}(N) \ll_{\gamma}M$.
		
		\item Si $N\leq M$ et s'il existe un épimorphisme $M/N\rightarrow M$,
		alors $N\ll_{\gamma} M$.
		
		\item Si $M/N$ est non nul et non cosingulier pour tout $N\leq M$ et si $f$
		est un endomorphisme surjectif de $M$ alors $f(N)\neq M$.

	\end{enumerate}
\end{thm}


L'exemple suivant montre que la classe des modules Hopfiens est une
sous-classe propre des modules $\gamma$-Hopfiens.

\begin{ex}\label{ex:4g}\cite{Alg}
	
	Soit $G=\mathbb{Z}_{p^{\infty}}$. Puisque $G$ est creux alors tout
	sous-groupe propre de $G$ est $\gamma$-superflu, donc $G$ est un groupe
	$\gamma$-Hopfien. Mais $G$ n'est pas Hopfien car la multiplication par
	$p$ induit un épimorphisme de $G$ qui n'est pas un isomorphisme.
\end{ex}

\begin{lem}\label{lem:abde}\cite{Alg}
	Soit $M$ un $A$-module et soit $N$ un sous-module de $M$, alors les
	assertions suivantes sont équivalentes:
	\begin{enumerate}
		
		\item $N\ll_{\gamma} M$.
		
		\item Si $X+N=M$, alors $X$ est un facteur direct de $M$ avec $M/X$ est un module semi simple
		cosingulier.
		
	\end{enumerate}
	
\end{lem}


\begin{thm}\label{thm:g}\cite{Alg} Soit $M$ un $A$-module, alors les
	assertions suivantes sont équivalentes:
	\begin{enumerate}
		
		\item $M$ est $\gamma$-Hopfien,
		
		\item Pour tout module $X$ tel qu'il existe un épimorphisme $M\rightarrow M\oplus X$, alors $X$ est semi simple cosingulier.

	\end{enumerate}
\end{thm}




\begin{thm}\label{thm:ggg}\cite{Alg} Soit $A$ un anneau. Alors les assertions suivantes sont équivalentes:
	\begin{enumerate}
		
		\item Tout $A$-module est $\gamma$-Hopfien,
		
		\item Tout $A$-module projectif est $\gamma$-Hopfien,
		
		\item Tout $A$-module libre est $\gamma$-Hopfien,
		
		\item $A$ est semi simple cosingulier.
		
	\end{enumerate}
\end{thm}

Il est clair que chaque module Hopfien généralisé est
$\gamma$-Hopfien. L'exemple suivant montre que la réciproque n'est pas vraie
en général. De plus, cet exemple montre aussi qu'un module
$\gamma$-Hopfien n'est pas en général Dedekind Fini.

\begin{ex}\label{ex:exg}\cite{Alg}
	
	Soit $A$ un anneau semi simple cosingulier. Alors d'après le
	théorème~\ref{thm:ggg}, $M=A^{(\mathbb{N})}$ est un $A$-module
	$\gamma$-Hopfien. Puisque $A^{(\mathbb{N})}\cong
	A^{(\mathbb{N})}\oplus A^{(\mathbb{N})}$ et $A^{(\mathbb{N})}\neq
	0$, $M$ n'est pas un module Hopfien généralisé (Dedekind
	Fini) (voir \cite[Corollaire 1.4]{Gh}).
\end{ex}

\begin{thm}\label{thm:morg}\cite{Alg}
	La propriété $\gamma$-Hopfien est préservée par l'équivalence de Morita.
\end{thm}

\begin{cor}\label{cor:mmg}\cite{Alg}
	Soit $n\geq 2$. Alors les assertions suivantes sont équivalentes
	pour un anneau $A$:
	
	(1) Tout $A$-module engendré par $n$ éléments est
	$\gamma$-Hopfien.
	
	(2) Tout $M_{n}(A)$-module cyclique est $\gamma$-Hopfien.
	
\end{cor}

\begin{thm}\cite{Alg} Soit $M$ un $A$-module, alors les
	assertions suivantes sont équivalentes:
	\begin{enumerate}
		
		\item $M$ est $\gamma$-Hopfien.
		
		\item Il existe un sous-module $\gamma$-superflu totalement invariant $N$ de $M$ tel que $M/N$ est $\gamma$-Hopfien.

	\end{enumerate}
\end{thm}


\begin{cor}\cite{Alg}
	Soit $M$ un module faiblement co-Hopfien. Si $M$ vérifie CCA sur
	les sous-modules non $\gamma$-superflus $N$ tel que $M/N$
	est faiblement co-Hopfien, alors $M$ est $\gamma$-Hopfien.
\end{cor}

\begin{prop}
	Soit $M$ un $A$-module. Si $M$ vérifie CCD sur les sous-modules
	non $\gamma$-superflus alors il est $\gamma$-Hopfien.
	
\end{prop}

\begin{prop}\cite{Alg}
	Soit $M$ un $A$-module vérifie la propriété suivante, pour
	tout endomorphisme $f$ de $M$ il existe un entier $n\geq 1$ tel que
	$Kerf^{n}\cap Imf^{n}\ll_{\gamma}M$. Alors $M$ est $\gamma$-Hopfien.
	
\end{prop}



Dans le corollaire suivant, on donne une caractérisation d'un anneau $A$ dans lequel tout $A$-module libre de type fini est $\gamma$-Hopfien.

\begin{cor}\cite{Alg}
	Soit $A$ un anneau. Alors les assertions suivantes
	sont équivalentes:
	\begin{enumerate}
		
		\item Tout $A$-module libre de type fini est $\gamma$-Hopfien.
		
		\item Tout $A$-module projectif de type fini est $\gamma$-Hopfien.
		
		\item $M_{n}(A)$ est un  $M_{n}(A)$-module $\gamma$-Hopfien pour tout
		$n\geq 1$.
		
	\end{enumerate}
	
\end{cor}


\begin{prop}\cite{Alg}
	Soit $M$ un module semi Hopfien, si $M$ est co-Hopfien alors il
	est $\gamma$-Hopfien.
\end{prop}

\section{Modules Jacobson Hopfiens}

\begin{defn}\cite{Ali}
	
	Un $A$-module $M$ est appelé Jacobson Hopfien, si pour tout
	endomorphisme surjectif $f$ de $M$, le noyau de $f$ est
	Jacobson-superflu dans $M$ $(Ker f \ll_{J} M)$.
	
\end{defn}
Le résultat suivant présente une caractérisation des
Modules Jacobson Hopfiens.

\begin{thm}\label{thm:ja}\cite{Ali} Soit $M$ un $A$-module, alors les
	assertions suivantes sont équivalentes:
	\begin{enumerate}
		
		\item $M$ est Jacobson Hopfien.
		
		\item Pour tout endomorphisme surjectif $f$ de $M$, si $N\ll_{J} M$,
		alors $f^{-1}(N) \ll_{J}M$.
		
		\item Pour tout épimorphisme $f : M/N\rightarrow M$, on a $N\ll_{J} M$.
		
		\item Si $M/N$ est non nul et $\mathrm{Rad}(M/N)=M/N$ pour tout $N\leq M$ et si $f$
		est un endomorphisme surjectif de $M$ alors $f(N)\neq M$.

	\end{enumerate}
\end{thm}


L'exemple suivant montre que la classe des modules Hopfiens est une
sous-classe propre des modules Jacobson Hopfiens.
\begin{ex}\label{ex:4}\cite{Ali}
	
	Soit $M=\mathbb{Z}_{p^{\infty}}$. Puisque tout sous-module de $M$ est Jacobson-superflu dans $M$ car $M$ est creux, alors il est clair que $M$ est Jacobson Hopfien mais $M$ n'est pas Hopfien. Noter que la multiplication par $p$ induit un épimorphisme de $G$ qui n'est pas un isomorphisme.
	
\end{ex}

\begin{rem}\label{rem:rem}\cite{Ali}
	Selon les définitions, tout module creux est Jacobson Hopfien, mais la réciproque n'est pas vraie en général. Noter que $M=\mathbb{Z}_{6}$ est un $\mathbb{Z}$-module semi simple n'est pas creux. Puisque pour tout module semi simple $M$ on a $\mathrm{Rad}(M)=0$, alors tout sous-module propre est Jacobson-superflu dans $M$ mais $M$ n'a aucun sous-module superflu non nul.
	
\end{rem}
\begin{lem}\label{lem:abderr}\cite{Ali}
	Soit $M$ un $A$-module et soit $K$ un sous-module de $M$, alors les
	assertions suivantes sont équivalentes:
	\begin{enumerate}
		
		\item $K\ll_{J} M$.
		
		\item Si $X+K=M$, alors $X$ est un facteur direct de $M$ avec $M/X$ est un module semi simple.
		
	\end{enumerate}
	
\end{lem}


\begin{thm}\label{thm:jaja}\cite{Ali} Soit $M$ un $A$-module, alors les
	assertions suivantes sont équivalentes:
	\begin{enumerate}
		
		\item $M$ est Jacobson Hopfien.
		
		\item Pour tout module $X$ tel qu'il existe un épimorphisme $M\rightarrow M\oplus X$, alors $X$ est semi simple.

	\end{enumerate}
\end{thm}



\begin{thm}\label{thm:jajaja}\cite{Ali} Soit $M$ un module (quasi-)projectif et $f\in \mathrm{End}(M)$.
	Alors les assertions suivantes sont équivalentes:
	\begin{enumerate}
		
		\item $M$ est Jacobson Hopfien.
		
		\item Si $f$ est un épimorphisme, alors $Ker(f)$ est semi simple.

	\end{enumerate}
\end{thm}


\begin{thm}\label{thm:jajajaja}\cite{Ali} Soit $A$ un anneau. Alors les assertions suivantes sont équivalentes:
	\begin{enumerate}
		
		\item Tout $A$-module est Jacobson Hopfien.
		
		\item Tout $A$-module projectif est Jacobson Hopfien.
		
		\item Tout $A$-module libre est Jacobson Hopfien.
		
		\item $A$ est semi simple.
		
	\end{enumerate}
\end{thm}

Il est clair que chaque module Hopfien généralisé est
Jacobson Hopfien. L'exemple suivant montre que la réciproque n'est pas vraie
en général. De plus, cet exemple montre aussi qu'un module
Jacobson Hopfien n'est pas en général Dedekind Fini.

\begin{ex}\label{ex:ex}\cite{Ali}
	
	Soit $A$ un anneau semi simple. Alors d'après le
	théorème~\ref{thm:jajajaja}, $M=A^{(\mathbb{N})}$ est un $A$-module
	Jacobson Hopfien. Puisque $A^{(\mathbb{N})}\cong
	A^{(\mathbb{N})}\oplus A^{(\mathbb{N})}$ et $A^{(\mathbb{N})}\neq
	0$, $M$ n'est pas Hopfien généralisé (Dedekind
	Fini) (voir \cite[Corollaire 1.4]{Gh}).
\end{ex}


\begin{thm}\cite{Ali} Soit $M$ un $A$-module, alors les
	assertions suivantes sont équivalentes:
	\begin{enumerate}
		
		\item $M$ est Jacobson Hopfien.
		
		\item Il existe un sous-module Jacobson-superflu totalement invariant $N$ de $M$ tel que $M/N$ est Jacobson Hopfien.

	\end{enumerate}
\end{thm}

\begin{prop}\label{prop:44}\cite{Ali}
	Soit $N$  un sous-module totalement invariant de $M$ tel que $M/N$
	est Hopfien. Si $N$ est Jacobson Hopfien alors $M$ est Jacobson Hopfien.
\end{prop}

\begin{prop}\label{pr:5}\cite{Ali}
	Soit $M$ un $A$-module. Si $M$ vérifie CCA sur les sous modules
	non Jacobson-superflus alors il est Jacobson Hopfien.
\end{prop}

\begin{prop}\cite{Ali}
	Soit $M$ un $A$-module. Si $M$ vérifie CCD sur les sous modules
	non Jacobson-superflus alors il est Jacobson Hopfien.
	
\end{prop}

\begin{prop}\cite{Ali}
	Soit $M$ un $A$-module vérifie la propriété suivante, pour
	tout endomorphisme $f$ de $M$ il existe un entier $n\geq 1$ tel que
	$Kerf^{n}\cap Imf^{n}\ll_{J}M$. Alors $M$ est Jacobson Hopfien.
	
\end{prop}


\begin{exs}
	\begin{enumerate}
		
		\item Tout sous-module propre de module semi simple $M$ est Jacobson-superflu, alors pour
		tout endomorphisme $f$ de $M$ il existe un entier $n\geq 1$ tel que
		$Kerf^{n}\cap Imf^{n}\ll_{J}M$. Donc $M$ est Jacobson Hopfien.
		
		\item Si $M$ est noethérien, alors pour tout endomorphisme
		$f$ de $M$ il existe un entier $n\geq 1$ tel que $Kerf^{n}\cap Imf^{n}=0$. Donc $M$ est Jacobson Hopfien.
		
	\end{enumerate}
	
\end{exs}

Dans le corollaire suivant, on donne une caractérisation d'un anneau $A$ dans lesquels tout $A$-module libre de type fini est Jacobson Hopfien.

\begin{cor}\cite{Ali}
	Soit $A$ un anneau. Alors les assertions suivantes
	sont équivalentes:
	\begin{enumerate}
		
		\item Tout $A$-module libre de type fini est Jacobson Hopfien.
		
		\item Tout $A$-module projectif de type fini est Jacobson Hopfien.
		
		\item $M_{n}(A)$ est un  $M_{n}(A)$-module Jacobson Hopfien pour tout
		$n\geq 1$.
		
	\end{enumerate}
	
\end{cor}

\begin{prop}
	Soit $M$ un module semi Hopfien, si $M$ est co-Hopfien alors il
	est Jacobson Hopfien.
\end{prop}


\section{Propriétés des extensions polynomiales}
Soit $M$ un $A$-module. D'après \cite{Va}, nous allons rappeler brièvement les
définitions des modules $M[x]$ et $M[x]/(x^{n+1})$. Les éléments de $M[x]$ sont des sommes formelles
de la forme $a_{0} + a_{1}x +...+ a_{k}x^{k}$ avec $k$ un entier
supérieur ou égal à $0$ et $a_{i }\in M$. On note cette
somme par $\sum^{k}_{i=1}a_{i}x^{i}$. Pour l'addition en ajoutant
les coefficients correspondants. La structure de $A[x]$-module est
définie par
\begin{center}
	$(\sum^{k}_{i=0}\lambda_{i}x^{i}).(\sum^{z}_{j=0}a_{j}x^{j}) =
	\sum^{k+z} _{\mu=0}c_{\mu}x^{\mu},$
\end{center}
où $c_{\mu} = \sum_{i+j=\mu}\lambda_{i}a_{j}$, pour tout
$\lambda_{i }\in A$, $a_{j }\in M$.

Tout élément non nul $\beta$ de $ M[x]$ s'écrit
uniquement sous la forme $(\sum^{l}_{i=k}m_{i}x^{i})$ avec $l \geq k
\geq 0$, $m_{i}\in M$, $m_{k}\neq 0$ et $m_{l}\neq 0$. Dans ce cas,
nous nous référons à $k$ comme l'ordre de $\beta$, $l$
comme le degré de $\beta$, et $m_{k}$ comme le coefficient initial
de $\beta$.

Soit $n$ un entier supérieur ou égal à $0$ et
\begin{center}
	$I_{n+1} = \{0\}\cup \{\beta ; 0\neq \beta\in A[x]$, l'ordre de
	$\beta \geq n+1\}$.
\end{center}
Alors $I_{n+1}$ est un idéal bilatéral de $A[x]$. L'anneau
quotient $A[x]/I_{n+1}$ est appelé l'anneau polynomial
tronqué, tronqué au degré $n+1$. Si $A$ est unitaire, $I_{n+1}$ est l'idéal
engendré par $x^{n+1}$. Si $A$ n'est pas unitaire, nous désignerons "symboliquement" l'anneau
$A[x]/I_{n+1}$ par $A[x]/(x^{n+1})$. Tout élément de
$A[x]/(x^{n+1})$ s'écrit uniquement sous la forme
$(\sum^{n}_{i=0}\lambda_{i}x^{i})$ avec $\lambda_{i }\in A$.

Soit
\begin{center}
	$D_{n+1} = \{0\}\cup \{\beta ; 0\neq \beta\in M[x]$, l'ordre de
	$\beta \geq n+1\}$.
\end{center}
Alors $D_{n+1}$ est un  $A[x]$-sous-module de $M[x]$. Puisque
$I_{n+1}M[x]\subset D_{n+1}$, on dit que $A[x]/(x^{n+1})$ agit sur
$M[x]/D_{n+1}$. On note le module $M[x]/D_{n+1}$ par
$M[x]/(x^{n+1})$. L'action de $A[x]/(x^{n+1})$ sur $M[x]/(x^{n+1})$
est donnée par
\begin{center}
	$(\sum^{n}_{i=0}\lambda_{i}x^{i}).(\sum^{n}_{j=0}a_{j}x^{j}) =
	\sum^{n} _{\mu=0}c_{\mu}x^{\mu},$
\end{center}
où $c_{\mu} = \sum_{i+j=\mu}\lambda_{i}a_{j}$, pour tout
$\lambda_{i }\in A$, $a_{j }\in M$.

Tout élément non nul $\beta$ de $ M[x]/D_{n+1}$ s'écrit
uniquement sous la forme $(\sum^{n}_{i=k}m_{i}x^{i})$ avec $n \geq k
\geq 0$, $m_{i}\in M$, $m_{k}\neq 0$. Dans ce cas, nous nous
référons $k$ comme l'ordre de $\beta$ et $m_{k}$ comme le
coefficient initial de $\beta$.

Le $A[x_{1},...,x_{k}]/(x_{1}^{n_{1}+1},...,x_{k}^{n_{k}+1})$-module
$M[x_{1},...,x_{k}]/(x_{1}^{n_{1}+1},...,x_{k}^{n_{k}+1})$ est
défini de la même manière.

\begin{lem}\label{lem:sss}\cite[Lemme 2.1]{Ga}.
	Soit $M$ un $A$-module et $K \ll M$. Alors $K[x]/(x^{n+1}) \ll
	M[x]/(x^{n+1})$ comme $A[x]/(x^{n+1})$-modules, où $n \geq 0$.
\end{lem}


\begin{lem}\label{lem:ssss}\cite[Lemme 1.7]{Va4}.
	Soit $N$ un $A$-sous-module de $M$. Alors les assertions suivantes
	sont équivalentes:
	
	\begin{enumerate}
		
		\item $N$ est un $A$-module essentiel dans $M$.
		
		\item $N[x]$ est un $A[x]$-module essentiel dans $M[x]$.
		
		\item $N[x]/(x^{n+1})$ est un $A[x]/(x^{n+1})$-module essentiel dans $M[x]/(x^{n+1})$.
		
	\end{enumerate}
\end{lem}
\begin{thm}\cite{Pal}
	Soit $M$ un $A$-module. Si $M[x]/(x^{n+1})$ est un
	$A[x]/(x^{n+1})$-module faiblement Hopfien, alors $M$ est un
	$A$-module faiblement Hopfien.
\end{thm}



\begin{thm}\cite{Pal}
	Soit $M$ un $A$-module. Si $M[x]/(x^{n+1})$ est un
	$A[x]/(x^{n+1})$-module co-Hopfien généralisé, alors $M$
	est un  $A$-module co-Hopfien généralisé.
\end{thm}

\begin{thm}\cite{Pal}
	Soit $M$ un $A$-module. Si
	$M[x_{1},...,x_{k}]/(x_{1}^{n_{1}+1},...,x_{k}^{n_{k}+1})$ est un
	$A[x_{1},...,x_{k}]/(x_{1}^{n_{1}+1},...,x_{k}^{n_{k}+1})$-module
	faiblement Hopfien (resp, co-Hopfien généralisé), alors
	$M$ est un $A$-module faiblement Hopfien (resp, co-Hopfien
	généralisé).
\end{thm}

\begin{lem}\cite{Jor}
	Soit $M$ un $A$-module et soit $N$ un sous-module de $M$. Si $N[x]/(x^{n+1})$
	est un facteur direct de $M[x]/(x^{n+1})$, alors $N$ est un facteur
	direct de $M$.
\end{lem}

\begin{thm}\cite{Jor}
	Soit $M$ un $A$-module. Si $M[x]/(x^{n+1})$ est un
	$A[x]/(x^{n+1})$-module semi Hopfien, alors $M$ est un $A$-module  semi
	Hopfien.
\end{thm}

\begin{thm}\cite{Jor}
	Soit M un $A$-module. Si $M[x]/(x^{n+1})$ est un
	$A[x]/(x^{n+1})$-module semi co-Hopfien, alors M est un A-module
	semi co-Hopfien.
\end{thm}

\begin{thm}\cite{Jor}
	Soit M un $A$-module. Si
	$M[x_{1},...,x_{k}]/(x_{1}^{n_{1}+1},...,x_{k}^{n_{k}+1})$ est un
	$A[x_{1},...,x_{k}]/(x_{1}^{n_{1}+1},...,x_{k}^{n_{k}+1})$-module
	semi Hopfien (resp, semi co-Hopfien), alors M est un A-module semi
	Hopfien (resp, semi co-Hopfien).
\end{thm}

\begin{lem}\cite{Fer}
	Soient $M$ un $A$-module et $K \ll_{\mu} M$. Alors $K[x]/(x^{n+1})
	\ll_{\mu} M[x]/(x^{n+1})$ comme $A[x]/(x^{n+1})$-modules, où
	$n\geq0$.
	
\end{lem}


\begin{thm}\cite{Fer}
	Soit $M$ un $A$-module. Alors $M[x]/(x^{n+1})$ est un
	$A[x]/(x^{n+1})$-module $\mu$-Hopfien
	si et seulement si $M$ est un $A$-module $\mu$-Hopfien.
\end{thm}

\begin{cor}\cite{Fer}
	Soit $M$ un $A$-module. Alors
	$M[x_{1},...,x_{k}]/(x_{1}^{n_{1}+1},...,x_{k}^{n_{k}+1})$ est
	$A[x_{1},...,x_{k}]/(x_{1}^{n_{1}+1},...,x_{k}^{n_{k}+1})$-module
	$\mu$-Hopfien si et seulement si $M$ est un $A$-module
	$\mu$-Hopfien.
\end{cor}


\end{document}